\documentclass[reqno,12pt]{amsart}
\usepackage{amsmath,amsthm,amssymb,amsfonts,amscd}
\usepackage{epsfig}
\usepackage{color}
\usepackage[all]{xy}
\usepackage{tikz}
\theoremstyle{plain}
\newtheorem{theorem}{Theorem}

\newtheorem{proposition}[theorem]{Proposition}
\newtheorem*{theorem*}{Theorem}
\newtheorem*{conjecture*}{Conjecture}

\theoremstyle{definition}

\newtheorem{definition}[theorem]{Definition}

\usepackage{graphics, setspace}

\newcommand{\breakingcomma}{%
  \begingroup\lccode`~=`,
  \lowercase{\endgroup\expandafter\def\expandafter~\expandafter{~\penalty0 }}}
  
\usepackage{breqn}

\usepackage{amstext} 
\usepackage{array}
\newcolumntype{L}{>{$}l<{$}}

\allowdisplaybreaks

\newcommand{\CC}{{\mathbb{C}}}

\newcommand{\QQ}{{\mathbb{Q}}}

\newcommand{\ZZ}{{\mathbb{Z}}}

\newcommand{\SL}{{\mathrm{SL}}}
\newcommand{\epi}{{\bf e}}

\newcommand{\Jac}{\mathrm{Jac}}
\newcommand{\Fix}{\mathrm{Fix}}

\newcommand{\id}{\mathrm{id}}

\newcommand{\bx}{{\bf x}}
\newcommand{\hess}{\mathrm{hess}}
\newcommand{\age}{\mathrm{age}}

\newcommand{\QR}[2]{
\left.\raisebox{0.5ex}{\ensuremath{#1}}
\ensuremath{\mkern-3mu}\middle/\ensuremath{\mkern-3mu}
\raisebox{-0.5ex}{\ensuremath{#2}}\right.}

\def\A{{\mathcal A}}

\begin{document}
\title{Orbifold Jacobian algebras for exceptional unimodal singularities}
\date{\today}
\author{Alexey Basalaev}
\address{Universit\"at Mannheim, Lehrsthul f\"ur Mathematik VI, Seminargeb\"aude A 5, 6, 68131 Mannheim, Germany}
\email{basalaev@uni-mannheim.de}
\author{Atsushi Takahashi}
\address{Department of Mathematics, Graduate School of Science, Osaka University, 
Toyonaka Osaka, 560-0043, Japan}
\email{takahashi@math.sci.osaka-u.ac.jp}
\author{Elisabeth Werner}
\address{Leibniz Universit\"at Hannover, Welfengarten 1, 30167 Hannover, Germany}
\email{werner.elisabeth@math.uni-hannover.de }
\begin{abstract}
This note shows that the orbifold Jacobian algebra associated to each invertible polynomial defining an exceptional unimodal singularity 
is isomorphic to the (usual) Jacobian algebra of the Berglund--H\"{u}bsch transform of an invertible polynomial 
defining the strange dual singularity in the sense of Arnold. 
\end{abstract}
\maketitle

\section{Introduction}
Exceptional unimodal singularities consist of $14$ isolated hypersurface singularities --- $Q_{10}$, $Q_{11}$, $Q_{12}$, $S_{11}$, $S_{12}$, $U_{12}$, $Z_{11}$, $Z_{12}$, $Z_{13}$, $W_{12}$, $W_{13}$, $E_{12}$, $E_{13}$ and $E_{14}$ in the Arnold's notation (see \cite{AGV85}). 
Arnold observed a "strange duality" in this class of singularities, the Dolgachev numbers (a triple of algebraically defined positive integers ) of one singularity are equal to the Garbrielov numbers (a triple of positive integers associated to a Coxeter--Dynkin diagram) of another one and vice versa. It is now naturally understood as one of mirror symmetry phenomena (cf. \cite{ET11} and references therein).

Let the polynomial $f \in \CC[x_1,\dots,x_N]$ be {\it invertible} (see Definition~\ref{def:invertible}). For such polynomials $f$ one can associate a priori new polynomial $\tilde f \in \CC[x_1,\dots,x_N]$, that is also invertible, called {\it Berglund--H\"ubsch transpose} of $f$ (see Section~\ref{OJA} for details). 

For any two exceptional unimodal singularities that are strange dual by Arnold there is a particular choice of the polynomials $f_1,f_2$ representing them such that both polynomials are invertible and $f_1 = \widetilde{f_2}$.
This was first observed in \cite{KY95}, where the authors show the coincidence of the elliptic genera of dual pairs up to sign, 
and also plays an essential role in \cite{ET11} for a precise formulation and generalization of Arnold's strange duality.
However, the choice of an invertible polynomial, representing an exceptional unimodal singularity is not unique in general (we list all possible choices of an invertible polynomial, representing an exceptional unimodal singularity in Table~\ref{tab: exceptional unimodal invertible}).

For an invertible polynomial $f \in \CC[x_1,\dots,x_N]$ and its symmetry group $G_f^\SL$ (see Section~\ref{OJA}), let
$\Jac(f,G_f^\SL)$ stand for the orbifold Jacobian algebra of the pair $(f,G_f^\SL)$ and $\Jac(f)$ be the ``usual'' Jacobian (or local) algebra. 
We prove the following theorem.
\begin{theorem}\label{theorem: SD}
Let $f_1,f_2 \in \CC[x,y,z]$ be invertible polynomials defining exceptional unimodal singularities (full list is given in Table~\ref{tab: exceptional unimodal invertible}). 
There exists a Frobenius algebra isomorphism
\[
\Jac(f_1)=\Jac(f_1,\{\id\}) \cong \Jac(\widetilde f_2,G^\SL_{\widetilde f_2})
\]
if and only if the associated singularities of $f_1$ and $f_2$ are strange dual to each other in the sense of Arnold. 
Here $\widetilde f_2$ is the Berglund--H\"ubsch transpose of $f_2$.
\end{theorem}

For a fixed singularity and different choices of the invertible polynomial $f_2$ representing it, the function $\tilde f_2$ can have different symmetry groups $G_{\tilde f_2}^\SL$ and even different Milnor numbers. In particular for $U_{12}$ one will get the symmetry groups $\{\id\}$, $\ZZ/2\ZZ$, $\ZZ/3\ZZ$ and Milnor numbers $12$, $12$, $15$ by $\tilde f_2$. The algebra $\Jac(f_1)$ in Theorem~\ref{theorem: SD} will still be the same up to isomorphism. Hence Theorem~\ref{theorem: SD} shows many non--trivial isomorphism.

\begin{table}
\begin{center}
\begin{tabular}{c||c|c|c||c}
Type & f (v1) & f (v2) & f (v3) 
& Strange dual \\ \hline
$Q_{10}$ & $x^4+y^3+x z^2$ & --  & -- 
&$E_{14}$ \\
$Q_{11}$ & $x^3 y+y^3+x z^2$ & --  & -- 
&$Z_{13}$ \\
$Q_{12}$ & $x^3 z+y^3+x z^2$ & $x^5+y^3+x z^2$ & --  
&$Q_{12}$ \\
$S_{11}$ & $x^4+y^2 z+x z^2$ &  --& -- 
&$W_{13}$ \\
$S_{12}$ & $x^3 y+y^2 z+x z^2$ & -- & -- 
&$S_{12}$ \\
$U_{12}$ & $x^4+y^3+z^3$ & $x^4+y^3+z^2 y$& $x^4+y^2 z+z^2 y$ 
&$U_{12}$ \\
$Z_{11}$ & $x^5+x y^3+z^2$ &  -- & -- 
&$E_{13}$ \\
$Z_{12}$ & $y x^4+x y^3+z^2$ &  --  & -- 
&$Z_{12}$ \\
$Z_{13}$ & $x^3 z+x y^3+z^2$ & $x^6+y^3 x+z^2$  & -- 
&$Q_{11}$ \\
$W_{12}$ & $x^5+y^2 z+z^2$ & $x^5+y^4+z^2$  & -- 
&$W_{12}$ \\
$W_{13}$ & $y x^4+y^2 z+z^2$ & $x^4 y+y^4+z^2$ & --  
&$S_{11}$ \\
$E_{12}$ & $x^7+y^3+z^2$ & -- & -- 
&$E_{12}$ \\
$E_{13}$ & $y^3+y x^5+z^2$ &  --  & -- 
&$Z_{11}$ \\
$E_{14}$ & $x^4 z+y^3+z^2$ & $x^8+y^3+z^2$ &  -- 
&$Q_{10}$ \\ \hline
\end{tabular}
\caption{This table shows all possible invertible polynomials, representing an exceptional unimodal singularity of a given type.}
\label{tab: exceptional unimodal invertible}
\end{center} 
\end{table}

It is worth to mention that Theorem~\ref{theorem: SD} is compatible with the mirror symmetry. 
Let $\A^{\mathrm{FJRW}}(f_2, \langle  g_{f_2} \rangle)$ stand for the so-called FJRW ring, the analogue of quantum cohomology ring associated to the pair $(f_2, \langle  g_{f_2} \rangle)$ and $g_{f_2}=(\epi\left[w_x/d\right], \epi\left[w_y/d\right],\epi\left[w_z/d\right])$ be the so-called 
exponential grading operator where $w_x,w_y,w_z, d$ are weights of variables $x,y,z$ and the polynomial $f_2=f_2(x,y,z)$ (see Section~\ref{OJA} for the notation).
An isomorphism of Frobenius algebras $ \Jac(f_1,\{\id\}) \cong \A^{\mathrm{FJRW}}(f_2, \langle  g_{f_2} \rangle) $ is obtained in \cite{KP+}.
As a corollary to Theorem~\ref{theorem: SD}, we get
\[
\A^{\mathrm{FJRW}}(f_2, \langle  g_{f_2} \rangle) \cong \Jac(\tilde{f_2},G^\SL_{\tilde{f_2}}),
\]
which is expected classical mirror symmetry isomorphism.

It's important to note that similar results are obtained in an apparently different context, 
in the study of matrix factorizations \cite{CRR16} and \cite{NR16}. 
We expect that the Hochschild cohomology group of the category of $G$-equivariant matrix factorizations will naturally yield the relationship between theirs and ours.
We hope to elaborate on this subject in the near future.

\begin{sloppypar}

{\bf Acknowledgements}.\
The first named author is partially supported by the DGF grant He2287/4--1 (SISYPH).
The second named author is supported by JSPS KAKENHI Grant Number JP16H06337, JP26610008.
We are grateful to Wolfgang Ebeling for fruitful discussions.
\end{sloppypar}

\section{Orbifold Jacobian algebra of an invertible polynomial}\label{OJA}

\subsection{Invertible polynomials}\label{section: invertible polynomials}
For a non-negative integer $n$ and $f=f(x_1,\dots, x_n)\in\CC[x_1,\dots, x_n]$ a polynomial, 
the {\em Jacobian algebra} $\Jac(f)$ of $f$ is a $\CC$-algebra defined as 
\begin{equation}
\Jac(f)= \QR{\CC[x_1,\dots, x_n]}{\left(\frac{\partial f}{\partial x_1},\dots,\frac{\partial f}{\partial x_n}\right)}.
\end{equation}
If $\Jac(f)$ is a finite-dimensional $\CC$-algebra, then set $\mu_f:=\dim_\CC\Jac(f)$ and call it the {\em Milnor number} of $f$. 
In particular, if $n=0$ then $\Jac(f)=\CC$ and $\mu_f=1$.

The {\em Hessian} of $f$ is defined as 
\begin{equation}
\hess(f):=\det \left(\frac{\partial^2 f}{\partial x_{i} \partial x_{j}}\right)_{i,j=1,\dots,n}.
\end{equation}
In particular, if $n=0$ then $\hess(f)=1$.

A polynomial $f\in\CC[x_1,\dots, x_N]$ is called a {\em weighted homogeneous} polynomial  
if there are positive integers $w_1,\dots ,w_N$ and $d$ such that 
\begin{equation}
f(\lambda^{w_1} x_1, \dots, \lambda^{w_N} x_N) = \lambda^d f(x_1,\dots ,x_N)
\end{equation}
for all $\lambda \in \CC^\ast$.
A weighted homogeneous polynomial $f$ is called {\em non-degenerate}
if it has at most an isolated critical point at the origin in $\CC^N$, equivalently, if the Jacobian algebra $\Jac(f)$ of $f$ 
is finite-dimensional.
\begin{definition}\label{def:invertible}
A non-degenerate weighted homogeneous polynomial $f\in\CC[x_1,\dots, x_N]$ is called {\em invertible} if 
the following conditions are satisfied.
\begin{itemize}
\item 
The number of variables ($=N$) coincides with the number of monomials in the polynomial $f$, namely, 
\begin{equation}
f(x_1,\dots ,x_N)=\sum_{i=1}^N c_i\prod_{j=1}^Nx_j^{E_{ij}}
\end{equation}
for some coefficients $c_i\in\CC^\ast$ and non-negative integers 
$E_{ij}$ for $i,j=1,\dots, N$.
\item
The matrix $E:=(E_{ij})$ is invertible over $\QQ$.
\item
The polynomial $f$ and the {\em Berglund--H\"{u}bsch transpose} $\widetilde f$ of $f$ defined by
\begin{equation}
\widetilde{f}(x_1,\dots ,x_N):=\sum_{i=1}^N c_i\prod_{j=1}^N x_j^{E_{ji}}
\end{equation}
are non-degenerate.
\end{itemize}
\end{definition}
\begin{definition}
The {\em group of maximal diagonal symmetries} of an invertible polynomial $f(x_1,\dots, x_N)$ is defined as 
\begin{equation}
G_f := \left\{ (\lambda_1, \ldots , \lambda_N )\in(\CC^\ast)^N \, \left| \, f(\lambda_1x_1, \ldots, \lambda_N x_N) = 
f(x_1, \ldots, x_N) \right\} \right..
\end{equation}
We shall always identify $G_f$ with the subgroup of diagonal matrices of ${\rm GL}(N;\CC)$.
Set 
\begin{equation}
G_f^{\rm SL} := G_f\cap {\rm SL}(N;\CC).
\end{equation}
Each element $g\in G_f$ has a unique expression of the form
\begin{equation}\label{G-notation}
g={\rm diag}\left(\epi\left[\frac{a_1}{r}\right], \dots,\epi\left[\frac{a_N}{r}\right]\right)
\quad \mbox{with } 0 \leq a_i < r,
\end{equation}
where $\epi\left[\alpha\right] := \exp(2 \pi \sqrt{-1} \alpha)$ and $r$ is the order of $g$. We use the notation $(a_1/r,\dots , a_N/r)$ for the element $g$.
The {\em age} of $g$ is defined as the rational number 
\begin{equation}
{\rm age}(g) := \frac{1}{r}\sum_{i=1}^N a_i. 
\end{equation}
Note that the ${\rm age}(g)$ is an integer if $g\in G_f^\SL$.
\end{definition}

\subsection{Orbifold Jacobian algebra}
Let $f=f(x_1,\dots, x_N)$ be an invertible polynomial and $G$ a subgroup of $G_f^\SL$.  
{\it A $G$--twisted Jacobian algebra of $f$} $\Jac'(f,G)$, which exists and is uniquely defined up to an isomorphism by  \cite[Theorem~21]{BTW16}, is given as follows.

As a $\CC$-vector space, $\Jac'(f,G)$ is given by
\begin{equation}\label{eq: JacFG as v.sp.}
\Jac'(f,G) = \bigoplus_{g \in G} \Jac(f^g)\widetilde v_g,
\end{equation}
where $f^g := f|_{\Fix(g)}$, $\Fix(g) \subseteq \CC^N$ is the fixed locus of $g$ and $\widetilde v_g$ is a generator (a formal letter) attached to each $g\in G$. Note that $f^g$ is also an invertible polynomial 
and there is a surjective map $\Jac(f)\longrightarrow \Jac(f^g)$ (\cite[Proposition~5]{ET13} and \cite[Proposition~7]{BTW16}).
It is also important that $\Jac'(f,G)$ is equipped with a $\ZZ/2\ZZ$-grading according to the parity of $N-N_g$, the codimension of the fixed locus $\Fix(g)$, for each $g\in G$.

We are now ready to introduce the product structure on $\Jac'(f,G)$. 
For simplicity, we assume that $G$ is a cyclic group whose order is a prime number.

For each pair $(g,h)$ of elements in $G$ and $\phi(\bx), \psi(\bx)\in \Jac(f)$, it is defined as follows:
\begin{itemize}
\item 
Suppose that $\Fix(g)\cup \Fix(h) \cup \Fix(gh) = \CC^N$.
Then
\begin{equation}\label{prod str}
[\phi(\bx)]\widetilde v_g\circ [\psi(\bx)]\widetilde v_h 
:=(-1)^{\frac{1}{2}(N-N_{g})(N-N_{g}-1)}\cdot \epi\left[-\frac{1}{2}\age(g)\right]\cdot  [\phi(\bx)\psi(\bx) H_{g,h}]\widetilde v_{gh},
\end{equation}
where $H_{g,h}\in\CC[x_1,\dots, x_N]$ is defined by the following equation in $\Jac(f^{gh})$
\begin{equation}\label{H and hessians}
\frac{1}{\mu_{f^{g\cap h}}}[\hess(f^{g\cap h})H_{g,h}]=\frac{1}{\mu_{f^{gh}}}[\hess(f^{gh})],
\end{equation}
and here $f^{g\cap h}$ is an invertible polynomial given by the restriction $f|_{\Fix(g)\cap\Fix(h)}$ of $f$ to 
the locus $\Fix(g)\cap \Fix(h)$.
\item 
Suppose that $\Fix(g)\cup \Fix(h) \cup \Fix(gh) \ne \CC^N$.
Then $[\phi(\bx)]\widetilde v_g\circ [\psi(\bx)]\widetilde v_h :=0$.
\end{itemize}
This completes the definition of $\Jac'(f,G)$. 
It is easy to see that $\widetilde v_\id=[1]\widetilde v_\id$ is the identity of $\Jac'(f,G)$.

Note that we have a natural action of $G$ on $\Jac(f^g)$ for any $g \in G$ and 
that the product structure is invariant under the $G$-action.
\begin{definition}
Let $f$ and $G$ be as above. The $G$-invariant $\ZZ/2\ZZ$-graded subalgebra $\Jac(f,G) := \left( \Jac'(f,G) \right)^G$ is called 
the {\em orbifold Jacobian algebra} of $(f,G)$. 
\end{definition}
An important property of this algebra is the following
\begin{proposition}[\cite{BTW16}]
The algebra $\Jac(f,G)$ is a $\ZZ/2\ZZ$-graded commutative Frobenius algebra. 
Namely, there is an even non-degenerate pairing $\eta_{f,G}$ such that
\begin{align}
& \eta_{f,G} \left( X \circ Y, Z \right)  = \eta_{f,G} \left( X, Y \circ Z \right),\quad X,Y,Z \in \Jac'(f,G),\\
& \eta_{f,G} \left( \widetilde v_\id, [\hess(f)]\widetilde v_\id \right)  = |G|\cdot \mu_f.\label{residue}
\end{align}
\end{proposition}

\section{Proof of Theorem~\ref{theorem: SD}}
The proof of Theorem~\ref{theorem: SD} is done by direct calculation. In what follows let the notation be as in Theorem~\ref{theorem: SD}.

Skipping the trivial cases when $f_1=\widetilde f_2$ and $G_{\widetilde f_2}^\SL=\{\id\}$, to prove the theorem we only need to show that $\Jac(f_1) \cong \Jac(\widetilde f_2,G_{\widetilde f_2}^\SL)$ for each row of Table~\ref{tab: full list} on page \pageref{tab: full list}.

\begin{table}  
\begin{tabular}{L|L||L|L|L}
\text{Type of $f_1$} & f_1 & \widetilde f_2 & G_{\widetilde f_2}^\SL &  \text{Type of $f_2$} \\ \hline
\mathrm{E_{14}} &  x_1^8+x_2^3+x_3^2  & x_1^4 x_2+x_2^2+x_3^3 & \{\id\} & \mathrm{Q_{10}} 
\\
\mathrm{Q_{10}} & x_1^4+x_2^3+x_1 x_3^2 & x_1^8+x_2^3+x_3^2 & \langle(1/2,0,1/2)\rangle &\mathrm{E_{14}} 
\\
\mathrm{Q_{11}} & x_1^3 x_2+x_2^3+x_1 x_3^2  & x_1^6 x_2+x_2^3+x_3^2 & \langle(1/2,0,1/2)\rangle &\mathrm{Z_{13}}   
\\
\mathrm{Q_{12}} & x_2^3+x_1^3 x_3+x_1 x_3^2  & x_1^5 x_2+x_2^2+x_3^3 & \langle (1/2,1/2,0) \rangle &\mathrm{Q_{12}}
\\
\mathrm{Q_{12}} & x_1^5+x_2^3+x_1 x_3^2  & x_1^3+x_2^3 x_3+x_2 x_3^2 & \{\id\} &\mathrm{Q_{12}}
\\
\mathrm{Q_{12}} &x_1^5+x_2^3+x_1 x_3^2  & x_1^5 x_2+x_2^2+x_3^3 & \langle(1/2,1/2,0)\rangle &\mathrm{Q_{12}}  
\\
\mathrm{S_{11}} &x_1^4+x_2^2 x_3+x_1 x_3^2  & x_1^4+x_1 x_2^4+x_3^2 & \langle(0,1/2,1/2)\rangle &\mathrm{W_{13}}
\\
\mathrm{U_{12}} &x_1^4+x_2^3+x_3^3  & x_1^4+x_2^3+x_3^3& \langle (0,2/3,1/3)\rangle &\mathrm{U_{12}}
\\
\mathrm{U_{12}}  &x_1^4+x_2^3+x_3^3  & x_1^4+x_2^3 x_3+x_3^2& \langle  (0,1/2,1/2)\rangle &\mathrm{U_{12}}
  \\
\mathrm{U_{12}}  &x_1^4+x_2^3+x_3^3  & x_1^4+x_2^2 x_3+x_2 x_3^2 & \{\id\} &\mathrm{U_{12}}
  \\
\mathrm{U_{12}}  &x_1^4+x_2^3+x_2 x_3^2  & x_1^4+x_2^3+x_3^3& \langle (0,2/3,1/3)\rangle &\mathrm{U_{12}}
  \\
\mathrm{U_{12}}  &x_1^4+x_2^3+x_2 x_3^2  & x_1^4+x_2^3 x_3+x_3^2& \langle  (0,1/2,1/2)\rangle &\mathrm{U_{12}}
  \\
\mathrm{U_{12}}  &x_1^4+x_2^3+x_2 x_3^2  & x_1^4+x_2^2 x_3+x_2 x_3^2 & \{\id\} &\mathrm{U_{12}}
  \\
\mathrm{U_{12}}  &x_1^4+x_2^2 x_3+x_2 x_3^2  & x_1^4+x_2^3+x_3^3& \langle (0,2/3,1/3)\rangle &\mathrm{U_{12}}
  \\
\mathrm{U_{12}}  &x_1^4+x_2^2 x_3+x_2 x_3^2  & x_1^4+x_2^3 x_3+x_3^2& \langle  (0,1/2,1/2)\rangle &\mathrm{U_{12}}
\\
\mathrm{W_{12}}  &x_1^5+x_2^2 x_3+x_3^2  & x_1^5+x_2^4+x_3^2 & \langle  (0,1/2,1/2)\rangle &\mathrm{W_{12}}
  \\
\mathrm{W_{12}}  &x_1^5+x_2^4+x_3^2  & x_1^5+x_2^2+x_2 x_3^2 & \{\id\} &\mathrm{W_{12}}
  \\
\mathrm{W_{12}}  &x_1^5+x_2^4+x_3^2  & x_1^5+x_2^4+x_3^2& \langle  (0,1/2,1/2)\rangle &\mathrm{W_{12}}
\\
\mathrm{W_{13}}, &x_1^4 x_2+x_2^4+x_3^2  & x_1^4 x_2+x_2^2 x_3+x_3^2 & \{\id\} &\mathrm{S_{11}}
\\
\mathrm{Z_{13}}, & x_1^6+x_1 x_2^3+x_3^2  & x_1^3 x_2+x_2^2+x_1 x_3^3 & \{\id\} & \mathrm{Q_{11}}
\end{tabular}
\caption{To prove Theorem~\ref{theorem: SD}, we need to show $\Jac(f_1) \cong \Jac(\widetilde f_2, G_{\widetilde f_2}^\SL)$ for every row of this table.}
\label{tab: full list}
\end{table}

Further, note that if $G_{\widetilde f_2}^\SL = \{1\}$ and $f_1,\widetilde f_2$ do not coincide but belong to 
the same right-equivalence class, the proof follows since the Jacobian algebra is an invariant of the right-equivalence class. 
Therefore, it is enough to show the statement for each row of Table~\ref{tab: reduced list} 
on page \pageref{tab: reduced list}.
\begin{table}  
\begin{tabular}{L|L||L|L|L}
\text{Type of $f_1$} & f_1 & \widetilde f_2 & \text{Group } G_{\widetilde f_2}^\SL &  \text{Type of $f_2$} 
\\ \hline
\mathrm{Q_{10}} & x_1^4+x_2^3+x_1 x_3^2 & x_1^8+x_2^3+x_3^2 & \langle(1/2,0,1/2)\rangle &\mathrm{E_{14}} 
\\
\mathrm{Q_{11}} & x_1^3 x_2+x_2^3+x_1 x_3^2  & x_1^6 x_2+x_2^3+x_3^2 & \langle(1/2,0,1/2)\rangle &\mathrm{Z_{13}}   
\\
\mathrm{Q_{12}} & x_1^5+x_2^3+x_1 x_3^2  & x_1^5 x_2+x_2^2+x_3^3 & \langle (1/2,1/2,0) \rangle &\mathrm{Q_{12}}
\\
\mathrm{S_{11}} &x_1^4+x_2^2 x_3+x_1 x_3^2  & x_1^4+x_1 x_2^4+x_3^2 & \langle(0,1/2,1/2)\rangle &\mathrm{W_{13}}
\\
\mathrm{U_{12}} &x_1^4+x_2^3+x_3^3  & x_1^4+x_2^3+x_3^3& \langle (0,2/3,1/3)\rangle &\mathrm{U_{12}}
\\
\mathrm{U_{12}}  &x_1^4+x_2^3+x_2 x_3^2  & x_1^4+x_2^3 x_3+x_3^2& \langle  (0,1/2,1/2)\rangle &\mathrm{U_{12}}
  \\
\mathrm{W_{12}}  &x_1^5+x_2^4+x_3^2  & x_1^5+x_2^4+x_3^2& \langle  (0,1/2,1/2)\rangle &\mathrm{W_{12}}
\end{tabular}
\caption{It is enough to show $\Jac(f_1) \cong \Jac(\widetilde f_2, G_{\widetilde f_2}^\SL)$ for every row of this table to prove Theorem~\ref{theorem: SD}.}
\label{tab: reduced list}
\end{table}

\subsection{Computations}
From now on, we shall use the notation of Table~\ref{tab: full list} on page \pageref{tab: full list}.
In order to check that two algebras are isomorphic, we represent $\Jac(\widetilde f_2,G_{\widetilde f_2}^\SL)$ as 
a quotient algebra of a polynomial ring in three variables. 
Namely, we will compute relations in $\Jac(\widetilde f_2,G_{\widetilde f_2}^\SL)$ and show 
the existence of a surjective algebra homomorphism from $\Jac(f_1)$, which turns out to be an isomorphism due to 
the dimension reason.

\subsubsection{$Q_{10}$ and $E_{14}$}
For $\widetilde f_2 = x_1^8+x_2^3+x_3^2$, and $G_{\widetilde f_2}^\SL = \langle g \rangle =\langle(\frac{1}{2},0,\frac{1}{2})\rangle$, 
$\Jac(\widetilde f_2,G_{\widetilde f_2}^\SL)$ is a $10$-dimensional $\CC$-vector space, whose basis can be chosen as 
\begin{equation}
 \widetilde v_\id  , [x_1^2]\widetilde v_\id  ,  [x_1^4]\widetilde v_\id  ,  [x_1^6]\widetilde v_\id  ,  [x_2]\widetilde v_\id  ,  [x_1^2 x_2]\widetilde v_\id  ,  [x_1^4x_2]\widetilde v_\id  ,  [x_1^6 x_2]\widetilde v_\id  , \quad  \widetilde v_g ,  [x_2]\widetilde v_g.
\end{equation}
The only non-trivial non-zero products in $\Jac(\widetilde f_2,G_{\widetilde f_2}^\SL)$, calculated by \eqref{prod str}, are given by 
\begin{eqnarray*}
& & [x_1^2]\widetilde v_{\id}\circ [x_1^2]\widetilde v_{\id}= [x_1^4] \widetilde v_{\id}, \ 
[x_1^2]\widetilde v_{\id}\circ [x_1^4]\widetilde v_{\id} = [x_1^6 ]\widetilde v_{\id}, \ 
[x_1^2]\widetilde v_{\id}\circ [x_2]\widetilde v_{\id} = [x_1^2 x_2] \widetilde v_{\id} ,\\
& & [x_1^4] \widetilde v_{\id}\circ [x_2] \widetilde v_{\id} = [x_1^4 x_2] \widetilde v_{\id},\ 
[x_1^6] \widetilde v_{\id}\circ [x_2] \widetilde v_{\id} =  [x_1^6 x_2]\widetilde v_{\id},\ 
[x_2] \widetilde v_{\id}\circ \widetilde v_{g} =  [x_2]\widetilde v_{g},\\
& & [x_1^2]\widetilde v_{\id}\circ [x_1^2 x_2] \widetilde v_{\id} =  [x_1^4 x_2] \widetilde v_{\id}, \ 
[x_1^4] \widetilde v_{\id}\circ [x_1^2 x_2] \widetilde v_{\id} = [x_1^6 x_2] \widetilde v_{\id},\\
& &  [x_1^2]\widetilde v_{\id}\circ [x_1^4 x_2] \widetilde v_{\id} =  [x_1^6 x_2] \widetilde v_{\id}, \ 
\widetilde v_{g}^2 = 16 [x_1^6] \widetilde v_{\id},\ 
\widetilde v_{g}\circ [x_2] \widetilde v_{g} = 16 [x_1^6 x_2] \widetilde v_{\id},
\end{eqnarray*}
which show that $[x_1^2]\widetilde v_\id$, $[x_2]\widetilde v_\id$, $\widetilde v_g$ generate $\Jac(\widetilde f_2,G_{\widetilde f_2}^\SL)$ and are subject to the following relations
\begin{equation}
16([x_1^2]\widetilde v_\id)^{\circ 3}-\widetilde v_g^{\circ 2}=0,\quad ([x_2]\widetilde v_\id)^{\circ 2}=0,\quad 
[x_1^2]\widetilde v_\id\circ \widetilde v_g=0.
\end{equation}

On the other hand, the Jacobian algebra $\Jac(f_1)$ is given by 
\begin{equation}
\Jac(f_1)= \QR {\CC[y_1,y_2,y_3]}{\left(4y_1^3+y_3^2,y_2^2,y_1y_3\right)}.
\end{equation}
Therefore, we have an algebra isomorphism 
\begin{equation}
\Jac(f_1)\stackrel{\cong}{\longrightarrow}\Jac(\widetilde f_2,G_{\widetilde f_2}^\SL),\quad 
y_1 \mapsto [x_1^2]\widetilde v_\id, \ y_2 \mapsto [x_2]\widetilde v_\id, \ y_3 \mapsto \frac{1}{2}\sqrt{-1}\widetilde v_g,
\end{equation}
which is, moreover, an isomorphism of Frobenius algebras since by~\eqref{residue} we have
\begin{align}
&\eta_{\widetilde f_2,G_{\widetilde f_2}^\SL} \left( \widetilde v_\id, [x_1^6x_2]\widetilde v_\id \right)  = 
\frac{1}{672}\cdot \eta_{\widetilde f_2,G_{\widetilde f_2}^\SL} \left( \widetilde v_\id, [\hess(\widetilde f_2)]\widetilde v_\id \right)  = \frac{2\cdot 14}{672}=\frac{1}{24},\\
& \eta_{f_1,\{\id\}} \left( \widetilde v_\id, [y_1^3y_2]\widetilde v_\id \right)  = 
\frac{1}{240}\cdot \eta_{f_1,\{\id\}} \left( \widetilde v_\id, [\hess(f_1)]\widetilde v_\id \right)  =\frac{1\cdot 10}{240} =\frac{1}{24}.
\end{align}

\subsubsection{$Q_{11}$ and $Z_{13}$}
For  $\widetilde f_2 = x_1^6 x_2+x_2^3+x_3^2$, and \( G^\SL_{\widetilde f_2} = \langle g \rangle = \langle( \frac{1}{2},0,\frac{1}{2})\rangle\), 
$\Jac(\widetilde f_2,G_{\widetilde f_2}^\SL)$ is a $11$-dimensional $\CC$-vector space, whose basis can be chosen as 
\begin{equation}
   \widetilde v_\id, [x_1^2]\widetilde v_\id, [x_1^4]\widetilde v_\id, [x_2]\widetilde v_\id, [x_1^2x_2]\widetilde v_\id, [x_1^4x_2]\widetilde v_\id, [x_2^2]\widetilde v_\id, [x_1^2 x_2^2]\widetilde v_\id, [x_1^4x_2^2]\widetilde v_\id, \quad \widetilde v_g, [x_2]\widetilde v_g.
\end{equation}
The only non-trivial non-zero products in $\Jac(\widetilde f_2,G_{\widetilde f_2}^\SL)$, calculated by \eqref{prod str}, are given by
\begin{eqnarray*}
& & [x_1^2]\widetilde v_{\id}\circ [x_1^2]\widetilde v_{\id} = [x_1^4]\widetilde v_{\id}, [x_1^2]\widetilde v_{\id} \circ [x_1^4]\widetilde v_{\id} = -3 [x_2^2]\widetilde v_{\id} , [x_1^4]\widetilde v_{\id}\circ [x_1^4]\widetilde v_{\id} = -3 [x_1^2 x_2^2]\widetilde v_{\id} , \\
& & [x_1^2]\widetilde v_{\id} \circ [x_2] \widetilde v_{\id} = [x_1^2 x_2] \widetilde v_{\id}, \
[x_1^4]\widetilde v_{\id}\circ [x_2]\widetilde  v_{\id}) = [x_1^4 x_2]\widetilde v_{\id} ,  \
[x_2] \widetilde v_{\id}\circ [x_2] \widetilde v_{\id} = [x_2^2]\widetilde v_{\id} ,  \\
& & 
[x_2]\widetilde v_{\id} \circ \widetilde v_{g} = [x_2]\widetilde v_{g}, \
[x_1^2]\widetilde v_{\id} \circ [x_1^2 x_2]\widetilde  v_{\id}) = [x_1^4 x_2]\widetilde v_{\id}, \
[x_2]\widetilde v_{\id} \circ [x_1^2 x_2]\widetilde  v_{\id}) = [x_1^2 x_2^2]\widetilde v_{\id} \\
& & [x_1^2 x_2] \widetilde v_{\id}\circ [x_1^2 x_2] \widetilde v_{\id} = [x_1^4 x_2^2]\widetilde v_{\id}, \
[x_2]\widetilde  v_{\id} \circ [x_1^4 x_2]\widetilde  v_{\id} = [x_1^4 x_2^2]\widetilde v_{\id},\
[x_1^2]\widetilde v_{\id} \circ [x_2^2]\widetilde v_{\id} = [x_1^2 x_2^2]\widetilde v_{\id},\\
& & [x_1^4]\widetilde v_{\id} \circ [x_2^2]\widetilde v_{\id}) = [x_1^4 x_2^2]\widetilde v_{\id},\
[x_1^2]\widetilde v_{\id} \circ [x_1^2 x_2^2]\widetilde v_{\id}= [x_1^4 x_2^2]\widetilde v_{\id},\
\widetilde v_{g} \circ \widetilde v_{g} = 12 [x_1^4 x_2]\widetilde v_{\id} , \\
& & \widetilde v_{g} \circ [x_2]\widetilde v_{g} = 12 [x_1^4 x_2^2]\widetilde v_{\id}.
\end{eqnarray*}
which show that $[x_1^2]\widetilde v_\id$, $[x_2]\widetilde v_\id$, $\widetilde v_g$ generate $\Jac(\widetilde f_2,G_{\widetilde f_2}^\SL)$ and are subject to the following relations
\begin{align}
12([x_1^2]\widetilde v_\id)^{\circ 2} [x_2]\widetilde v_\id - (\widetilde v_g)^{\circ 2} = 0, \ ([x_1^2]\widetilde v_\id)^{\circ 3} + 3([x_2]\widetilde v_\id)^{\circ 2} = 0, \ [x_1^2]\widetilde v_\id \circ \widetilde v_g = 0.
\end{align}

On the other hand, the Jacobian algebra $\Jac(f_1)$ is given by 
\[
\Jac(f)=\QR {\CC[y_1,y_2,y_3]}{\left(3y_1^2y_2+y_3^2,y_1^3+3y_2^2,2y_1y_3\right)}
\]
Therefore, we have an algebra isomorphism 
\begin{equation}
\Jac(f_1)\stackrel{\cong}{\longrightarrow}\Jac(\widetilde f_2,G_{\widetilde f_2}^\SL),\quad 
y_1 \mapsto [x_1^2]\widetilde v_\id, \ y_2 \mapsto [x_2]\widetilde v_\id, \ y_3 \mapsto \frac{1}{2}\sqrt{-1}\widetilde v_g,
\end{equation}
which is, moreover, an isomorphism of Frobenius algebras since by~\eqref{residue} we have
\begin{align}
&\eta_{\widetilde f_2,G_{\widetilde f_2}^\SL} \left( \widetilde v_\id, [x_1^4x_2^2]\widetilde v_\id \right)  = \frac{1}{18},
\\
& \eta_{f_1,\{\id\}} \left( \widetilde v_\id, [y_1^2y_2^2]\widetilde v_\id \right)  = 
\frac{1}{198}\cdot \eta_{f_1,\{\id\}} \left( \widetilde v_\id, [\hess(f_1)]\widetilde v_\id \right) = \frac{1\cdot 11}{198} =\frac{1}{18}.
\end{align}

\subsubsection{$Q_{12}$ and $Q_{12}$}
For $\widetilde f_2 = x_1^5 x_2+x_2^2+x_3^3$, and $G_{\widetilde f_2}^\SL = \langle g \rangle =\langle(\frac{1}{2},\frac{1}{2},0)\rangle$, 
$\Jac(\widetilde f_2,G_{\widetilde f_2}^\SL)$ is a $12$-dimensional $\CC$-vector space, whose basis can be chosen as 
\begin{align}
  & \widetilde v_\id, [x_1^2] \widetilde v_\id, [x_1^4]\widetilde v_\id, [x_1x_2]\widetilde v_\id, [x_1^3 x_2]\widetilde v_\id, [x_3]\widetilde v_\id, [x_1^2 x_3]\widetilde v_\id,
\\
  & [x_1^4 x_3]\widetilde v_\id, [x_1 x_2 x_3]\widetilde v_\id, [x_1^3 x_2 x_3]\widetilde v_\id, \ \widetilde v_g, [x_3]\widetilde v_g.
\end{align}
The only non-trivial non-zero products in $\Jac(\widetilde f_2,G_{\widetilde f_2}^\SL)$, calculated by \eqref{prod str}, are given by 
\begin{eqnarray*}
& & [x_1^2]\widetilde v_{\id} \circ [x_1^2]\widetilde v_{\id} = [x_1^4]\widetilde v_{\id}, \
[x_1^2] \widetilde v_{\id} \circ [x_1^4]\widetilde v_{\id} = -2 [x_1 x_2]\widetilde v_{\id}, \
[x_1^4\widetilde v_{\id}]\circ x_1^4\widetilde v_{\id} = -2 [x_1^3 x_2]\widetilde v_{\id}, \\
& & 
[x_1^2]\widetilde v_{\id} \circ [x_1 x_2]\widetilde v_{\id} = [x_1^3 x_2]\widetilde v_{\id} ,\
[x_1^2]\widetilde v_{\id} \circ [x_3]\widetilde v_{\id} = [x_1^2 x_3]\widetilde v_{\id} , \
[x_1^4]\widetilde v_{\id} \circ [x_3] \widetilde v_{\id} = [x_1^4 x_3]\widetilde v_{\id} ,\\
& & 
[x_1x_2]\widetilde v_{\id} \circ  [x_3]\widetilde v_{\id} = [x_1x_2 x_3]\widetilde v_{\id} , \
[x_1^3 x_2]\widetilde v_{\id} \circ  [x_3]\widetilde v_{\id} = [x_1^3x_2 x_3]\widetilde v_{\id} ,\
[x_3]\widetilde v_{\id} \circ \widetilde v_{g} = [x_3]\widetilde v_{g} , \\
& &
[x_1^2]\widetilde v_{\id} \circ [x_1^2 x_3]\widetilde v_{\id} = [x_1^4 x_3]\widetilde v_{\id} ,\
[x_1^4]\widetilde v_{\id}  \circ [x_1^2 x_3]\widetilde v_{\id} = -2 [x_1 x_2 x_3]\widetilde v_{\id} ,\\
& & 
[x_1x_2]\widetilde v_{\id} \circ [x_1^2 x_3]\widetilde v_{\id} = [x_1^3 x_2 x_3]\widetilde v_{\id} ,\
[x_1^2]\widetilde v_{\id} \circ [x_1^4 x_3]\widetilde v_{\id} = -2 [x_1 x_2 x_3]\widetilde v_{\id} ,\\
& &
[x_1^4]\widetilde v_{\id} \circ [x_1^4 x_3]\widetilde v_{\id} = -2 [x_1^3x_2 x_3]\widetilde v_{\id},\ 
[x_1^2]\widetilde v_{\id} \circ [x_1 x_2 x_3]\widetilde v_{\id} = [x_1^3 x_2 x_3]\widetilde v_{\id},\\
& &
\widetilde v_{g} \circ \widetilde v_{g} = 10 [x_1^3 x_2]\widetilde v_{\id} ,\
\widetilde v_{g} \circ [x_3]\widetilde v_{g} = 10 [x_1^3 x_2 x_3]\widetilde v_{\id}.
\end{eqnarray*}
which show that $[x_1^2]\widetilde v_\id$, $[x_3]\widetilde v_\id$, $\widetilde v_g$ generate $\Jac(\widetilde f_2,G_{\widetilde f_2}^\SL)$ and are subject to the following relations
\begin{equation}
[x_1^2]\widetilde v_\id \circ \widetilde v_g = 0, \ (\widetilde v_g)^{\circ 2} - 5([x_1^2]\widetilde v_\id)^{\circ 4} = 0,\ ([x_3]\widetilde v_{\id})^{\circ 2} = 0
\end{equation}

On the other hand, the Jacobian algebra $\Jac(f_1)$ is given by 
\begin{equation}
\Jac(f_1)= \QR {\CC[y_1,y_2,y_3]}{\left(5y_1^4+y_3^2,3y_2^2,2y_1y_3\right)}.
\end{equation}
Therefore, we have an algebra isomorphism 
\begin{equation}
\Jac(f_1)\stackrel{\cong}{\longrightarrow}\Jac(\widetilde f_2,G_{\widetilde f_2}^\SL),\quad 
y_1 \mapsto [x_1^2]\widetilde v_\id, \
y_2 \mapsto \frac{1}{4} [x_3]\widetilde v_\id, \
y_3 \mapsto \sqrt{-1}\widetilde v_g.
\end{equation}
which is, moreover, an isomorphism of Frobenius algebras since by~\eqref{residue} we have
\begin{align}
&\eta_{\widetilde f_2,G_{\widetilde f_2}^\SL} \left( \widetilde v_\id, [x_1^3x_2x_3]\widetilde v_\id \right) =\frac{1}{15},\\
- \frac{4}{10} \cdot & \eta_{f_1,\{\id\}} \left( \widetilde v_\id, [y_2y_3^2]\widetilde v_\id \right)  = 
\frac{4}{720}\cdot \eta_{f_1,\{\id\}} \left( \widetilde v_\id, [\hess(f_1)]\widetilde v_\id \right) = \frac{4\cdot 12}{720} =\frac{1}{15}.
\end{align}

\subsubsection{$S_{11}$ and $W_{13}$}
For $\widetilde f_2 = x_1^4+x_1 x_2^4+x_3^2$, and $G_{\widetilde f_2}^\SL = \langle g \rangle =\langle(0,\frac{1}{2},\frac{1}{2})\rangle$, 
$\Jac(\widetilde f_2,G_{\widetilde f_2}^\SL)$ is a $11$-dimensional $\CC$-vector space, whose basis can be chosen as 
\begin{align}
 \widetilde v_\id, [x_1]\widetilde v_\id, [x_1^2]\widetilde v_\id, [x_1^3]\widetilde v_\id, [x_2^2]\widetilde v_\id, [x_1x_2^2]\widetilde v_\id, [x_1^2 x_2^2]\widetilde v_\id, [x_1^3 x_2^2]\widetilde v_\id, \quad \widetilde v_g, [x_1]\widetilde v_g, [x_1^2]\widetilde v_g.
\end{align}
The only non-trivial non-zero products in $\Jac(\widetilde f_2,G_{\widetilde f_2}^\SL)$, calculated by \eqref{prod str}, are given by 
\begin{eqnarray*}
& & 
[x_1]\widetilde v_{\id} \circ [x_1]\widetilde v_{\id} = [x_1^2]\widetilde v_{\id}, \
[x_1]\widetilde v_{\id} \circ \widetilde v_{g} = [x_1]\widetilde v_{g},\
[x_1]\widetilde v_{\id} \circ [x_1]\widetilde v_{g} = [x_1^2]\widetilde v_{g},\\
& &
[x_1]\widetilde v_{\id} \circ [x_1^2]\widetilde v_{\id} = [x_1^3]\widetilde v_{\id} ,\
[x_1^2]\widetilde v_{\id} \circ \widetilde v_{g} = [x_1^2]\widetilde v_{g} ,\
[x_1]\widetilde v_{\id} \circ [x_2^2]\widetilde v_{\id} = [x_1 x_2^2]\widetilde v_{\id},\\
& &
[x_1^2]\widetilde v_{\id} \circ [x_2^2]\widetilde v_{\id} = [x_1^2 x_2^2]\widetilde v_{\id}, \
[x_1^3]\widetilde v_{\id} \circ [x_2^2]\widetilde v_{\id} = [x_1^3 x_2^2]\widetilde v_{\id} ,\
[x_2^2]\widetilde v_{\id} \circ [x_2^2]\widetilde v_{\id} = -4 [x_1^3]\widetilde v_{\id} ,\\
& & 
[x_1]\widetilde v_{\id} \circ [x_1 x_2^2]\widetilde v_{\id} = [x_1^2 x_2^2]\widetilde v_{\id} ,\
[x_1^2]\widetilde v_{\id} \circ [x_1 x_2^2]\widetilde v_{\id} = [x_1^3 x_2^2]\widetilde v_{\id} ,\\
& &
[x_1]\widetilde v_{\id} \circ [x_1^2 x_2^2]\widetilde v_{\id} = [x_1^3 x_2^2]\widetilde v_{\id} ,\
\widetilde v_{g} \circ \widetilde v_{g} = 8 [x_1 x_2^2]\widetilde v_{\id} ,\
\widetilde v_{g} \circ [x_1]\widetilde v_{g} = 8 [x_1^2 x_2^2]\widetilde v_{\id} ,\\
& &
[x_1]\widetilde v_{g} \circ [x_1]\widetilde v_{g} = 8 [x_1^3 x_2^2]\widetilde v_{\id}, \
\widetilde v_{g}\circ [x_1^2]\widetilde v_{g} = 8 [x_1^3 x_2^2]\widetilde v_{\id}.
\end{eqnarray*}
which show that $[x_1]\widetilde v_\id$, $[x_2^2]\widetilde v_\id$, $\widetilde v_g$ generate $\Jac(\widetilde f_2,G_{\widetilde f_2}^\SL)$ and are subject to the following relations
\begin{equation}
([x_2^2]\widetilde v_\id)^{\circ 2} + 4([x_1]\widetilde v_\id)^{\circ 3} = 0, \ [x_2^2]\widetilde v_\id \circ \widetilde v_g =0, \ (\widetilde v_g)^{\circ 2} - 8 [x_1]\widetilde v_\id \circ [x_2^2]\widetilde v_\id = 0.
\end{equation}

On the other hand, the Jacobian algebra $\Jac(f_1)$ is given by 
\begin{equation}
\Jac(f_1)= \QR {\CC[y_1,y_2,y_3]}{\left(4y_1^3+y_3^2,2y_2y_3,y_2^2+2y_1y_3\right)}.
\end{equation}
Therefore, we have an algebra isomorphism 
\begin{equation}
\Jac(f_1)\stackrel{\cong}{\longrightarrow}\Jac(\widetilde f_2,G_{\widetilde f_2}^\SL),\quad 
y_1\mapsto [x_1]\widetilde v_\id,\
y_2\mapsto \frac{1}{2}\sqrt{-1}\widetilde v_g,\
y_3\mapsto [x_2^2]\widetilde v_\id.
\end{equation}
which is, moreover, an isomorphism of Frobenius algebras since by~\eqref{residue} we have
\begin{align}
&\eta_{\widetilde f_2,G_{\widetilde f_2}^\SL} \left( \widetilde v_\id, [x_1^3x_2^2]\widetilde v_\id \right)  =\frac{1}{16},
\\
& \eta_{f_1,\{\id\}} \left( \widetilde v_\id, [y_1^3y_3]\widetilde v_\id \right)  = 
\frac{1}{176}\cdot \eta_{f_1,\{\id\}} \left( \widetilde v_\id, [\hess(f_1)]\widetilde v_\id \right)  = \frac{1 \cdot 11}{176} = \frac{1}{16}.
\end{align}

\subsubsection{$U_{12}$ and $U_{12}$, part 1}
For $\widetilde f_2 = x_1^4+x_2^3+x_3^3$, and $G_{\widetilde f_2}^\SL = \langle g \rangle =\langle(0,\frac{2}{3},\frac{1}{3})\rangle$, 
$\Jac(\widetilde f_2,G_{\widetilde f_2}^\SL)$ is a $12$-dimensional $\CC$-vector space, whose basis can be chosen as 
\begin{align}
& \widetilde v_\id, [x_1]\widetilde v_\id, [x_1^2]\widetilde v_\id, [x_2x_3]\widetilde v_\id, [x_1 x_2 x_3]\widetilde v_\id, [x_1^2 x_2 x_3]\widetilde v_\id, 
\\
& \widetilde v_{g^2}, [x_1]\widetilde v_{g^2}, [x_1^2]\widetilde v_{g^2}, \ \widetilde v_g, [x_1]\widetilde v_g, [x_1^2]\widetilde v_g.
\end{align}
The only non-trivial non-zero products in $\Jac(\widetilde f_2,G_{\widetilde f_2}^\SL)$, calculated by \eqref{prod str}, are given by 
\begin{eqnarray*}
& & 
[x_1]\widetilde v_{\id} \circ [x_1]\widetilde v_{\id} = [x_1^2]\widetilde v_{\id},\
[x_1]\widetilde v_{\id} \circ \widetilde v_{g^2} = [x_1]\widetilde v_{g^2} ,\
[x_1]\widetilde v_{\id} \circ [x_1]\widetilde v_{g^2} = [x_1^2]\widetilde v_{g^2},\\
& & 
[x_1]\widetilde v_{\id} \circ \widetilde v_{g} = [x_1]\widetilde v_{g},\
[x_1]\widetilde v_{\id} \circ [x_1]\widetilde v_{g} = [x_1^2]\widetilde v_{g} ,\
[x_1^2]\widetilde v_{\id} \circ \widetilde v_{g^2} = [x_1^2]\widetilde v_{g^2} ,\\
& &
[x_1^2\widetilde v_{\id}] \circ \widetilde v_{g} = [x_1^2]\widetilde v_{g} ,\
[x_1]\widetilde v_{\id} \circ [x_2 x_3]\widetilde v_{\id}) = [x_1 x_2 x_3]\widetilde v_{\id} ,\
[x_1^2]\widetilde v_{\id} \circ [x_2 x_3]\widetilde v_{\id} = [x_1^2 x_2 x_3]\widetilde v_{\id} ,\\
& &
[x_1]\widetilde v_{\id} \circ [x_1 x_2 x_3]\widetilde v_{\id}) = [x_1^2 x_2 x_3]\widetilde v_{\id} ,\
\widetilde v_{g} \circ \widetilde v_{g^2} = 9 [x_2 x_3]\widetilde v_{\id} ,\
\widetilde v_{g} \circ [x_1]\widetilde v_{g^2} = 9 [x_1 x_2 x_3]\widetilde v_{\id} ,\\
& &
\widetilde v_{g} \circ [x_1^2]\widetilde v_{g^2} = 9[x_1^2 x_2 x_3]\widetilde v_{\id} ,\
[x_1]\widetilde v_{g} \circ \widetilde v_{g^2} = 9 [x_1 x_2 x_3]\widetilde v_{\id} ,\
[x_1]\widetilde v_{g} \circ [x_1]\widetilde v_{g^2} = 9[x_1^2 x_2 x_3]\widetilde v_{\id} ,\\
& &
[x_1^2]\widetilde v_{g} \circ \widetilde v_{g^2} = 9 [x_1^2 x_2 x_3]\widetilde v_{\id}.
\end{eqnarray*}
which show that $[x_1]\widetilde v_\id$, $\widetilde v_g$, $\widetilde v_{g^2}$ generate $\Jac(\widetilde f_2,G_{\widetilde f_2}^\SL)$ and are subject to the following relations
\begin{equation}
([x_1]\widetilde v_\id)^{\circ 3} = 0, \ (\widetilde v_g)^{\circ 2} = 0, \ (\widetilde v_{g^2})^{\circ 2} = 0.
\end{equation}

On the other hand, the Jacobian algebra $\Jac(f_1)$ is given by 
\begin{equation}
\Jac(f_1)= \QR {\CC[y_1,y_2,y_3]}{\left(4y_1^3,3y_2^2,3y_3^2\right)}.
\end{equation}
Therefore, we have an algebra isomorphism 
\begin{equation}
\Jac(f_1)\stackrel{\cong}{\longrightarrow}\Jac(\widetilde f_2,G_{\widetilde f_2}^\SL),\quad 
y_1\mapsto \frac{1}{3} [x_1]\widetilde v_\id,\
y_2\mapsto \frac{1}{\sqrt{3}} \widetilde v_g,\
y_3\mapsto \frac{1}{\sqrt{3}} \widetilde v_{g^2}.
\end{equation}
which is, moreover, an isomorphism of Frobenius algebras since by~\eqref{residue} we have
\begin{align}
&\eta_{\widetilde f_2,G_{\widetilde f_2}^\SL} \left( \widetilde v_\id, [x_1^2x_2x_3]\widetilde v_\id \right)  =\frac{1}{12},\\
3 \cdot & \eta_{f_1,\{\id\}} \left( \widetilde v_\id, [y_1^2y_2y_3]\widetilde v_\id \right)  = 
\frac{3}{432}\cdot \eta_{f_1,\{\id\}} \left( \widetilde v_\id, [\hess(f_1)]\widetilde v_\id \right) =\frac{3\cdot 12}{432} =\frac{1}{12}.
\end{align}

\subsubsection{$U_{12}$ and $U_{12}$, part 2}
For $\widetilde f_2 = x_1^4+x_2^3 x_3+x_3^2$, and $G_{\widetilde f_2}^\SL = \langle g \rangle =\langle(0,\frac{1}{2},\frac{1}{2})\rangle$, 
$\Jac(\widetilde f_2,G_{\widetilde f_2}^\SL)$ is a $12$-dimensional $\CC$-vector space, whose basis can be chosen as 
\begin{align}
& \widetilde v_\id, [x_1]\widetilde v_\id, [x_1^2]\widetilde v_\id, [x_2^2]\widetilde v_\id, [x_1
x_2^2]\widetilde v_\id, [x_1^2 x_2^2]\widetilde v_\id, [x_2 x_3]\widetilde v_\id, [x_1 x_2 x_3]\widetilde v_\id, [x_1^2 x_2 x_3]\widetilde v_\id, 
\\
&\widetilde v_g, [x_1]\widetilde v_g, [x_1^2]\widetilde v_g.
\end{align}
The only non-trivial non-zero products in $\Jac(\widetilde f_2,G_{\widetilde f_2}^\SL)$, calculated by \eqref{prod str}, are given by 
\begin{eqnarray*}
& & 
[x_1]\widetilde v_{\id} \circ [x_1]\widetilde v_{\id} = [x_1^2]\widetilde v_{\id},\
[x_1]\widetilde v_{\id} \circ \widetilde v_{g} = [x_1]\widetilde v_{g} ,\
[x_1]\widetilde v_{\id} \circ [x_1]\widetilde v_{g} = [x_1^2]\widetilde v_{g},\\
& &
[x_1^2]\widetilde v_{\id} \circ \widetilde v_{g} = [x_1^2]\widetilde v_{g} ,\
[x_1]\widetilde v_{\id} \circ [x_2^2]\widetilde v_{\id} = [x_1 x_2^2]\widetilde v_{\id} ,\
[x_1^2]\widetilde v_{\id} \circ [x_2^2]\widetilde v_{\id} = [x_1^2 x_2^2]\widetilde v_{\id},\\
& &
[x_2^2]\widetilde v_{\id} \circ [x_2^2]\widetilde v_{\id} = -2 [x_2 x_3]\widetilde v_{\id} ,\
[x_1]\widetilde v_{\id} \circ [x_1 x_2^2]\widetilde v_{\id} = [x_1^2 x_2^2]\widetilde v_{\id} ,\\
& &
[x_2^2]\widetilde v_{\id} \circ  [x_1 x_2^2]\widetilde v_{\id} = -2[x_1 x_2x_3] \widetilde v_{\id} ,\
[x_1 x_2^2]\widetilde v_{\id} \circ [x_1 x_2^2]\widetilde v_{\id} = -2 [x_1^2 x_2 x_3]\widetilde v_{\id} ,\\
& &
[x_2^2] \widetilde v_{\id} \circ [x_1^2 x_2^2]\widetilde v_{\id} = -2 [x_1^2 x_2 x_3]\widetilde v_{\id} ,\
[x_1]\widetilde v_{\id} \circ [x_2 x_3]\widetilde v_{id} = [x_1 x_2 x_3]\widetilde v_{\id} ,\\
& & 
[x_1^2]\widetilde v_{\id} \circ [x_2 x_3]\widetilde v_{\id} = [x_1^2 x_2 x_3]\widetilde v_{\id} ,\
[x_1]\widetilde v_{\id} \circ [x_1 x_2 x_3]\widetilde v_{\id} = [x_1^2 x_2 x_3]\widetilde v_{\id} ,\\
& &
\widetilde v_{g} \circ \widetilde v_{g} = 6 [x_2 x_3]\widetilde v_{\id} , \
\widetilde v_{g} \circ [x_1]\widetilde v_{g} = 6[x_1 x_2 x_3] \widetilde v_{\id}, \\
& &
[x_1]\widetilde v_{g} \circ [x_1]\widetilde v_{g} = 6 [x_1^2 x_2 x_3]\widetilde v_{\id} ,\
\widetilde v_{g} \circ [x_1^2]\widetilde v_{g} = 6[x_1^2 x_2 x_3] \widetilde v_{\id}.
\end{eqnarray*}
which show that $[x_1]\widetilde v_\id$, $[x_2^2]\widetilde v_\id$, $\widetilde v_g$ generate $\Jac(\widetilde f_2,G_{\widetilde f_2}^\SL)$ and are subject to the following relations
\begin{equation}
([x_1]\widetilde v_\id)^{\circ 3} = 0, \ (\widetilde v_g)^{\circ 2} + 3([x_2^2]\widetilde v_\id)^{\circ 2}, \ [x_2^2]\widetilde v_\id \circ \widetilde v_g = 0.
\end{equation}

On the other hand, the Jacobian algebra $\Jac(f_1)$ is given by 
\begin{equation}
\Jac(f_1)= \QR {\CC[y_1,y_2,y_3]}{\left(4y_1^3,3y_2^2+y_3^2,2y_2y_3\right)}.
\end{equation}
Therefore, we have an algebra isomorphism 
\begin{equation}
\Jac(f_1)\stackrel{\cong}{\longrightarrow}\Jac(\widetilde f_2,G_{\widetilde f_2}^\SL),\quad 
y_1\mapsto \frac{1}{\sqrt{-6}} [x_1]\widetilde v_\id,\
y_2\mapsto [x_2^2]\widetilde v_\id,\
y_3\mapsto \widetilde v_{g}.
\end{equation}
which is, moreover, an isomorphism of Frobenius algebras since by~\eqref{residue} we have
\begin{align}
&\eta_{\widetilde f_2,G_{\widetilde f_2}^\SL} \left( \widetilde v_\id, [x_1^2x_2x_3]\widetilde v_\id \right)  =\frac{1}{12},
\\
\frac{6}{9} \cdot & \eta_{f_1,\{\id\}} \left( \widetilde v_\id, [y_1^2y_3^2]\widetilde v_\id \right)  = 
\frac{6}{9\cdot 96}\cdot \eta_{f_1,\{\id\}} \left( \widetilde v_\id, [\hess(f_1)]\widetilde v_\id \right) = \frac{6\cdot 12}{9\cdot 96} = \frac{1}{12}.
\end{align}

\subsubsection{$W_{12}$ and $W_{12}$}
For $\widetilde f_2 = x_1^5+x_2^4+x_3^2$, and $G_{\widetilde f_2}^\SL = \langle g \rangle =\langle(0,\frac{1}{2},\frac{1}{2})\rangle$, 
$\Jac(\widetilde f_2,G_{\widetilde f_2}^\SL)$ is a $12$-dimensional $\CC$-vector space, whose basis can be chosen as 
\begin{align}
\widetilde v_\id, [x_1]\widetilde v_\id, [x_1^2]\widetilde v_\id, [x_1^3]\widetilde v_\id, [x_2^2]\widetilde v_\id, [x_1 x_2^2]\widetilde v_\id, [x_1^2 x_2^2]\widetilde v_\id, [x_1^3 x_2^2]\widetilde v_\id, \quad \widetilde v_g, [x_1]\widetilde v_g, [x_1^2]\widetilde v_g, [x_1^3]\widetilde v_g.
\end{align}
The only non-trivial non-zero products in $\Jac(\widetilde f_2,G_{\widetilde f_2}^\SL)$, calculated by \eqref{prod str}, are given by 
\begin{eqnarray*}
& & 
[x_1]\widetilde v_{\id} \circ [x_1]\widetilde v_{\id} = [x_1^2]\widetilde v_{\id}, \
[x_1]\widetilde v_{\id} \circ \widetilde v_{g} = [x_1] \widetilde v_{g},\
[x_1]\widetilde v_{\id} \circ [x_1]\widetilde v_{g} = [x_1^2]\widetilde v_{g} ,\\
& &
[x_1] v_{\id} \circ [x_1^2]\widetilde v_{g} = [x_1^3]\widetilde v_{g},\
[x_1]\widetilde v_{\id} \circ [x_1^2]\widetilde v_{\id} = [x_1^3]\widetilde v_{\id},\
[x_1^2]\widetilde v_{\id} \circ \widetilde v_{g} = [x_1^2]\widetilde v_{g},\\
& &
[x_1^2]\widetilde v_{\id} \circ [x_1] \widetilde v_{g} = [x_1^3]\widetilde v_{g},\
[x_1^3]\widetilde v_{\id} \circ \widetilde v_{g} = [x_1^3]\widetilde v_{g}, \
[x_1]\widetilde v_{\id} \circ [x_2^2]\widetilde v_{\id} = [x_1 x_2^2]\widetilde v_{\id},\\
& & 
[x_1^2]\widetilde v_{\id} \circ [x_2^2]\widetilde v_{\id} =[x_1^2 x_2^2] \widetilde v_{\id} ,\
[x_1^3]\widetilde v_{\id} \circ [x_2^2]\widetilde v_{\id} = [x_1^3 x_2^2]\widetilde v_{\id} ,\
[x_1]\widetilde v_{\id} \circ [x_1 x_2^2]\widetilde v_{\id} = [x_1^2 x_2^2]\widetilde v_{\id},\\
& &
[x_1^2] \widetilde v_{\id} \circ [x_1 x_2^2]\widetilde v_{\id} = [x_1^3 x_2^2]\widetilde v_{\id} ,\
[x_1] \widetilde v_{\id} \circ [x_1^2 x_2^2] \widetilde v_{\id} = [x_1^3 x_2^2]\widetilde v_{\id} ,\ 
\widetilde v_{g} \circ \widetilde v_{g} = 8 [x_2^2]\widetilde v_{\id} , \\
& &
\widetilde v_{g} \circ [x_1]\widetilde v_g = 8[x_1 x_2^2] \widetilde v_{\id} , \
[x_1]\widetilde v_{g} \circ [x_1]\widetilde v_{g} = 8[x_1^2 x_2^2] \widetilde v_{\id} ,\
\widetilde v_{g} \circ [x_1^2]\widetilde v_{g} = 8 [x_1^2 x_2^2]\widetilde v_{\id},\\
& &
[x_1]\widetilde v_{g} \circ [x_1^2]\widetilde v_{g} = 8 [x_1^3 x_2^2]\widetilde v_{\id} ,\
\widetilde v_{g} \circ [x_1^3] \widetilde v_{g} = 8 [x_1^3 x_2^2]\widetilde v_{\id}.
\end{eqnarray*}
which show that $[x_1]\widetilde v_\id$ and $\widetilde v_g$ generate $\Jac(\widetilde f_2,G_{\widetilde f_2}^\SL)$ and are subject to the following relations
\begin{equation}
([x_1]\widetilde v_\id)^{\circ 4} = 0, \ (\widetilde v_{g})^{\circ 3} = 0
\end{equation}

On the other hand, the Jacobian algebra $\Jac(f_1)$ is given by 
\begin{equation}
\Jac(f_1)= \QR {\CC[y_1,y_2,y_3]}{\left(5y_1^4,4y_2^3,y_3\right)}.
\end{equation}
Therefore, we have an algebra isomorphism 
\begin{equation}
\Jac(f_1)\stackrel{\cong}{\longrightarrow}\Jac(\widetilde f_2,G_{\widetilde f_2}^\SL),\quad 
y_1\mapsto \frac{1}{2} [x_1]\widetilde v_\id,\
y_2\mapsto \frac{1}{\sqrt{2}}  v_g.
\end{equation}
which is, moreover, an isomorphism of Frobenius algebras since by~\eqref{residue} we have
\begin{align}
&\eta_{\widetilde f_2,G_{\widetilde f_2}^\SL} \left( \widetilde v_\id, [x_1^3x_2^2]\widetilde v_\id \right) =\frac{1}{20},
\\
\frac{16}{8} \cdot & \eta_{f_1,\{\id\}} \left( \widetilde v_\id, [y_1^3y_2^2]\widetilde v_\id \right)  = 
\frac{16}{8 \cdot 480}\cdot \eta_{f_1,\{\id\}} \left( \widetilde v_\id, [\hess(f_1)]\widetilde v_\id \right)  = \frac{16\cdot 12}{480} = \frac{1}{20}.
\end{align}

\subsection{Remark}
In order to visualize the statement of Theorem~\ref{theorem: SD} consider the following Figure~\ref{fig} on page~\pageref{fig}.
The nodes of this figure are the pairs $(f,G)$ where $f$ is an invertible polynomial and $G \subseteq G^\SL_f$. 
The edge between two nodes labeled by $(f_1,G_1)$ and $(f_2,G_2)$ respectively is drawn if and only if $\Jac(f_1,G_1) \cong \Jac(f_2,G_2)$.
All the pairs $(f,G)$ considered are those from Table~\ref{tab: full list}. 

\begin{figure}
  
\begin{tikzpicture}[yscale=-1]
\begin{scope}[every node/.style={rectangle,thick,draw}]
    \node (A) at (0,0) {$x^3 + y^4 + yz^2, \{\id\}$};
    \node (B) at (3,2) {$x^8+y^3+z^2,\ZZ/2\ZZ$};
    \node (C) at (-0.5,3.5) {$x^4 + xz^2 + y^3, \{\id\}$};
\end{scope}
\begin{scope}[
              every node/.style={fill=white,circle},
  	      every edge/.style={draw=black, thick}]
    \path [-] (A) edge (B);
    \path [-] (B) edge (C);
    \path [-] (A) edge (C);
\end{scope}

\begin{scope}[every node/.style={rectangle,thick,draw}]
    \node (D) at (8,0) {$x^4 + y^3 + z^3, \ZZ/3\ZZ$};
    \node (E) at (11,2) {$x^4+y^2z + y z^2,\{\id\}$};
    \node (F) at (7.5,3.5) {$x^4 + y^3z + z^2, \ZZ/2\ZZ$};
\end{scope}
\begin{scope}[
              every node/.style={fill=white,circle},
  	      every edge/.style={draw=black, thick}]
    \path [-] (F) edge (E);
    \path [-] (E) edge (D);
    \path [-] (F) edge (D);
\end{scope}
%
%
\begin{scope}[every node/.style={rectangle,thick,draw}]
    \node (G) at (0,6) {$x^5 + y^2 + yz^2, \{\id\}$};
    \node (H) at (3,8) {$x^5+y^4+z^2,\ZZ/2\ZZ$};
    \node (I) at (-0.5,9.5) {$x^5 + y^2z + z^2, \{\id\}$};
\end{scope}
\begin{scope}[
              every node/.style={fill=white,circle},
  	      every edge/.style={draw=black, thick}]
    \path [-] (G) edge (I);
    \path [-] (H) edge (G);
    \path [-] (I) edge (H);
\end{scope}
\begin{scope}[every node/.style={rectangle,thick,draw}]
    \node (J) at (8,6) {$x^4 + y^3 + z^3, \ZZ/3\ZZ$};
    \node (K) at (12,7.3) {$x^4+y^2z + y z^2,\{\id\}$};
    \node (L) at (7.5,9.5) {$x^4 + y^3z + z^2, \ZZ/2\ZZ$};
    \node (M) at (11.5,10.8) {$x^4 + y^3z + z^2, \ZZ/2\ZZ$};
\end{scope}
\begin{scope}[
              every node/.style={fill=white,circle},
  	      every edge/.style={draw=black, thick}]
    \path [-] (J) edge (K);
    \path [-] (L) edge (K);
    \path [-] (L) edge (J);
    \path [-] (M) edge (J);
    \path [-] (M) edge (K);
    \path [-] (M) edge (L);
\end{scope}
%
%
\begin{scope}[every node/.style={rectangle,thick,draw}]
    \node (O) at (0,12) {$x^5y + y^2 + z^3, \ZZ/2\ZZ$};
    \node (P) at (4,13.3) {$x^3 + y^3z + yz^2,\{\id\}$};
    \node (Q) at (-0.5,15.5) {$x^3z + xz^2 + y^3, \{\id\}$};
    \node (R) at (3.5,16.8) {$x^5 + xz^2 + y^3, \{\id\}$};
\end{scope}
\begin{scope}[
              every node/.style={fill=white,circle},
  	      every edge/.style={draw=black, thick}]
    \path [-] (O) edge (P);
    \path [-] (Q) edge (P);
    \path [-] (Q) edge (O);
    \path [-] (R) edge (O);
    \path [-] (R) edge (P);
    \path [-] (R) edge (Q);
\end{scope}
\begin{scope}[every node/.style={rectangle,thick,draw}]
    \node (S) at (8,12) {$x^3y + xz^3 + y^2, \{\id\}$};
    \node (T) at (12,13.3) {$x^3z+xy^3 + z^2,\{\id\}$};
    \node (U) at (8,14.5) {$x^6y + y^3 + z^2, \ZZ/2\ZZ$};
    \node (V) at (12,15.8) {$x^3y + xz^2 + y^3, \{\id\}$};;
    \node (W) at (8,17) {$x^4z + y^3 + z^2, \{\id\}$};
    \node (X) at (12,18.3) {$x^4y + y^2 + z^3, \{\id\}$};
\end{scope}
\begin{scope}[
              every node/.style={fill=white,circle},
  	      every edge/.style={draw=black, thick}]
    \path [-] (S) edge (T);
    \path [-] (U) edge (V);
    \path [-] (W) edge (X);
\end{scope}
\end{tikzpicture}
\caption{This figure depicts the isomorphisms between the different orbifold Jacobian algebras. The edge between two nodes labeled by $(f_1,G_1)$ and $(f_2,G_2)$ is drawn if and only if $\Jac(f_1,G_1) \cong \Jac(f_2,G_2)$. }
\label{fig}
\end{figure}
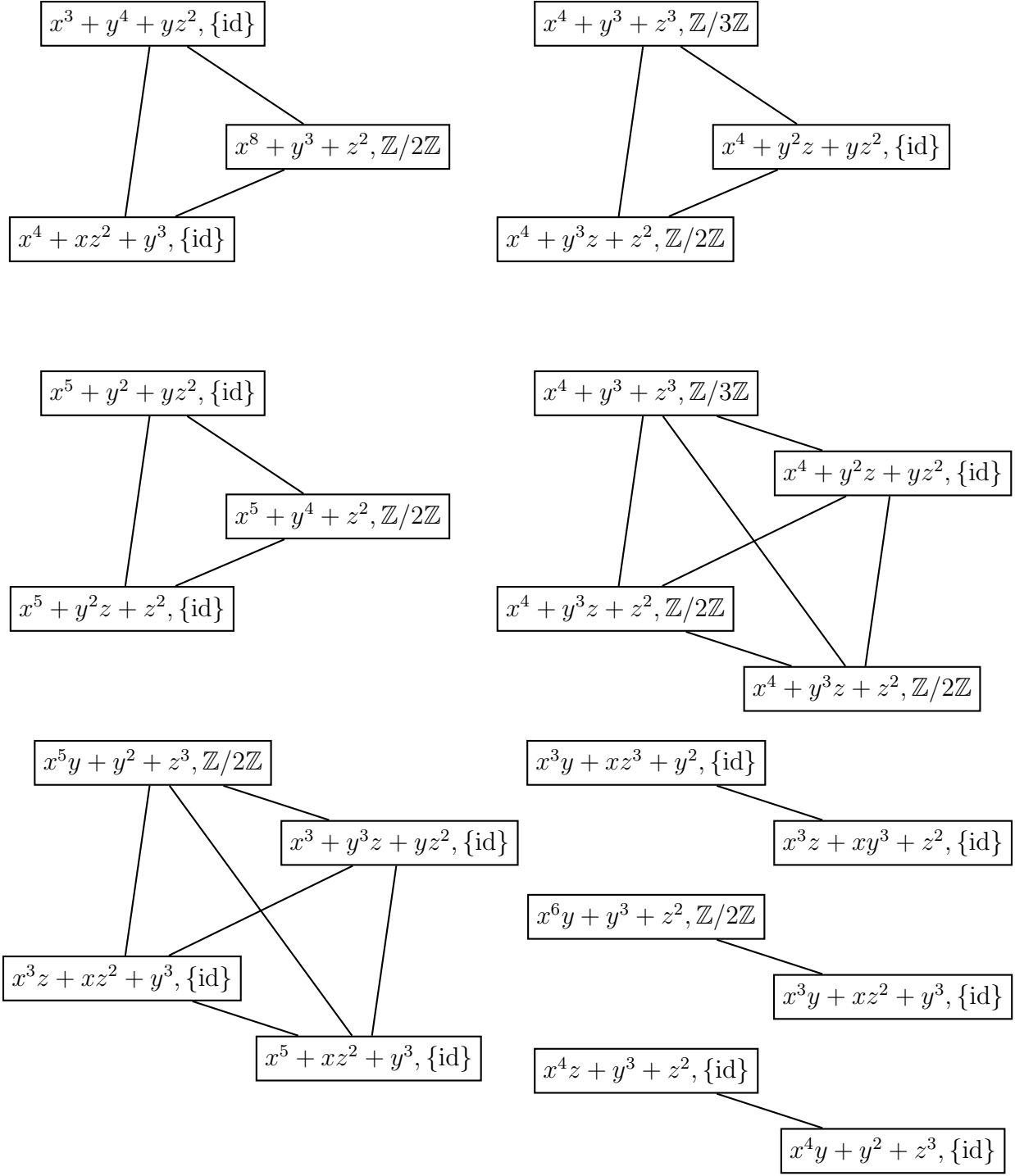


\end{document}